\documentclass[12pt]{article}
\usepackage[english]{babel}
\usepackage{amssymb, amsmath, amsfonts}
\usepackage{mathrsfs}
\newtheorem{Th}{Theorem}
\newtheorem{Lm}{Lemma}

\begin{document}
\begin{center}
{\bf ASYMPTOTIC BEHAVIOR OF LOCAL PARTICLES NUMBERS\\
IN BRANCHING RANDOM WALK}
\end{center}
\vskip0,5cm
\begin{center}
Ekaterina Vl. Bulinskaya\footnote{Lomonosov Moscow State
University.}$^,$\footnote{The work is partially supported by the
RFBR grant 10-01-00266.}
\end{center}
\vskip1cm

\begin{abstract}
Critical catalytic branching random walk on an integer lattice
$\mathbb{Z}^{d}$ is investigated for all $d\in\mathbb{N}$. The
branching may occur at the origin only and the start point is
arbitrary. The asymptotic behavior, as time grows to infinity, is
determined for the mean local particles numbers. The same problem is
solved for the probability of particles presence at a fixed lattice
point. Moreover, the Yaglom type limit theorem is established for
the local number of particles. Our analysis involves construction of
an auxiliary Bellman-Harris branching process with six types of
particles. The proofs employ the asymptotic properties of the
(improper) c.d.f. of hitting times with taboo. The latter notion was
recently introduced by the author for a non-branching random walk on
$\mathbb{Z}^{d}$.

\vskip0,5cm {\it Keywords and phrases}: critical branching random
walk, Bellman-Harris process with particles of six types, Yaglom
type conditional limit theorems, Kolmogorov's equations, random walk
on integer lattice, hitting time with taboo.

\vskip0,5cm {\it $2010$ AMS classification}: 60F05
\end{abstract}

\section{Introduction}

\emph{Catalytic branching random walk} (CBRW) on $d$-dimensional
integer lattice is a model of particles population evolution. We
recall its main features. Each particle independently of others may
perform random walk on $\mathbb{Z}^{d}$ and produce offsprings at
the source of branching located w.l.g. at the origin.
\emph{Symmetric branching random walk} (SBRW) on $\mathbb{Z}^{d}$
studied earlier, e.g., in \cite{ABY}, \cite{BY} and \cite{Y_UMN} is
a particular case of CBRW on $\mathbb{Z}^{d}$ (see \cite{Y_TVP}).

The model under consideration was proposed in \cite{VTY} for $d=1$
and studied for other $d\in\mathbb{N}$ in \cite{B_TVP},
\cite{B_TSP}, \cite{VT_Siberia} and \cite{Y_TVP}. The analysis of
CBRW in \cite{VT_Siberia} and \cite{Y_TVP} has shown that similarly
to many kinds of branching processes (see \cite{Sev}) CBRW on
$\mathbb{Z}^{d}$ is classified as supercritical, critical or
subcritical. According to \cite{Y_TVP}, the \emph{exponential
growth} (as time tends to infinity) of total number of particles in
population and local numbers of particles as well is characteristic
for the supercritical CBRW on $\mathbb{Z}^{d}$. The term local
refers to the (number of) particles located at a lattice point.

Quite different situation occurs for critical CBRW which is the main
object of study in this paper. For example, for $d=1$ or $d=2$ the
particles population degenerates with probability 1 but survives
with strictly positive probability for $d\geq3$ (see \cite{B_TSP},
\cite{MZ_Bulinskaya}, \cite{VT_Siberia} and \cite{VTY}). Moreover,
the total number of particles conditioned on non-degeneracy has
non-trivial discrete limit distribution, different for $d<3$ and
$d\geq3$ (see \cite{B_TVP}, \cite{B_TSP} and \cite{VTY}). Thus, in
the model of critical CBRW on $\mathbb{Z}^{d}$ the asymptotic
behavior (in time) of the total number of particles on the lattice
depends on dimension $d$ essentially and does not grow
exponentially. As for the local particles numbers in the model of
critical CBRW on $\mathbb{Z}^{d}$, earlier in
\cite{B_TVP}--\cite{LMJ}, \cite{VT_4}, \cite{VT_Siberia},
\cite{VT_DM} and \cite{VTY} there were only established the
asymptotic properties of the number of particles located at the
\emph{source of branching}. In particular, it turnes out that for
all $d\in\mathbb{N}$ the probability of the presence of particles at
the source of branching asymptotically vanishes. Notably, its
asymptotic behavior as well as limit laws for properly normalized
number of particles at the source of branching, conditioned on the
presence of particles at the origin, have different forms for
$d=1,2,3,4$ and $d\geq5$. Among the arising limit distributions one
can find exponential and discrete ones along with a mixture of such
laws.

In the model of critical CBRW on $\mathbb{Z}^{d}$ the behavior of
number of particles located at an \emph{arbitrary} point of the
lattice remained unknown. The present work completes the picture. We
study the asymptotic behavior in time of mean local particles
numbers and that of probability of particles presence at a fixed
point $y\neq{\bf 0}$ where ${{\bf
0}=(0,0,\ldots,0)\in\mathbb{Z}^{d}}$. All the more, we obtain a
conditional limit theorem for the properly normalized number of
particles at such point $y$. It should be emphasized that in
contrast to \cite{B_TVP}, \cite{VT_4}, \cite{VT_Siberia},
\cite{VT_DM} and \cite{VTY} we admit the start of CBRW at an
\emph{arbitrary} point $x\in\mathbb{Z}^{d}$ and not only at the
source of branching. Asymptotic properties of the number of
particles at ${\bf 0}$ for CBRW with an arbitrary start point were
investigated in \cite{LMJ}.

The structure of the rest of the paper is the following. In section
2 we describe the model in detail and formulate three main results.
Theorem \ref{T:m(t;x,y)} is proved in Section 3. Section 4 is
devoted to construction of the auxiliary Bellman-Harris branching
process. Thereupon we establish Theorems \ref{T:q(t;x,y)sim} and
\ref{T:limits} in Section 5.

\section{Main results}

Now we dwell on the definition of a critical CBRW on
$\mathbb{Z}^{d}$. At the initial time $t=0$ there is a single
particle on the lattice located at a point $x\in\mathbb{Z}^{d}$. If
$x\neq{\bf 0}$, the particle performs a continuous time random walk
until the time of the first hitting the origin. The random walk
outside the origin is symmetric, homogeneous, irreducible (i.e. a
particle passes from an arbitrary $u\in\mathbb{Z}^{d}$ to any
$\upsilon\in\mathbb{Z}^{d}$ with positive probability within a
finite time) and has a finite variance of jumps. Accordingly, we
assume this random walk be specified by an infinitesimal matrix
$A=(a(u,\upsilon))_{u,\;\upsilon\in\mathbb{Z}^{d}}$ such that
\begin{eqnarray*}
a(u,\upsilon)&=&a(\upsilon,u),\;\; a(u,\upsilon)=a({\bf
0},\upsilon-u):=a(\upsilon-u),\;\;
u,\upsilon\in\mathbb{Z}^{d},\\
\sum\nolimits_{\upsilon\in\mathbb{Z}^{d}}{a(\upsilon)}&=&0\;\;\mbox{where}\;\;
a({\bf 0})<0\;\;\mbox{and}\;\; a(\upsilon)\geq
0\;\;\mbox{if}\;\;\upsilon\neq{\bf
0},\;\;\sum\nolimits_{\upsilon\in\mathbb{Z}^{d}}{\|\upsilon\|^{2}a(\upsilon)}<\infty.
\end{eqnarray*}

If $x={\bf 0}$ or the particle has just hit the origin it spends
there an exponentially distributed time (with parameter 1).
Afterwards, it either dies with probability $\alpha\in(0,1)$
producing before the death a random number of offsprings $\xi$ or
leaves the source of branching with probability $1-\alpha$. In the
latter case the intensity of transition from the origin to a point
$\upsilon\neq{\bf 0}$ is given by
$$\overline{a}({\bf 0},\upsilon)=-(1-\alpha)\frac{a(\upsilon)}{a({\bf
0})}.$$
At the origin the branching is determined by a probability
generating function
$${f(s):={\sf
E}{s^{\xi}}=\sum\nolimits_{k=0}^{\infty}{f_{k}s^{k}}},\quad
s\in[0,1].$$ In \cite{VT_Siberia} CBRW on $\mathbb{Z}^{d}$ is called
\emph{critical} if the following relations hold
\begin{equation}\label{assumption_3}
\alpha
f'(1)+(1-\alpha)(1-h_{d})=1\quad\mbox{and}\quad\sigma^{2}:=\alpha\,f''(1)<\infty.
\end{equation}
Here $h_{d}$ is the probability of the event that a particle leaving
the origin will never return there. By the recurrence of a random
walk on $\mathbb{Z}$ and $\mathbb{Z}^{2}$ one has $h_{1}=h_{2}=0$.
It is well known that $h_{d}\in(0,1)$ for $d\geq3$.

Newborn particles are located at the origin at the birth moment.
They evolve according to the scheme described above independently of
each other as well as of the parent particles. The number of
particles located at a point $y\in\mathbb{Z}^{d}$ at time $t\geq0$
is denoted by $\mu(t;y)$.

The goal of the paper is three-fold. Firstly, we find the asymptotic
behavior (as $t\to\infty$) of the mean number of particles
$m(t;x,y):={\sf E}_{x}{\mu(t;y)}$ located at a point
$y\in\mathbb{Z}^{d}$, $y\neq{\bf 0}$, at time $t\geq0$ (everywhere
the index $x$ means that our CBRW starts at $x\in\mathbb{Z}^{d}$).
Secondly, we retrieve the asymptotic behavior of the probability
$q(t;x,y):={\sf P}_{x}(\mu(t;y)>0)$ of the presence of particles at
the point $y$ at time $t$. Thirdly, we establish a limit theorem for
properly normalized local numbers $\mu(t;y)$ conditioned on
$\mu(t;y)>0$ as $t\to\infty$.

To formulate the main results of the paper we introduce some more
notation. Let $p(t;x,y)$ be the transition probability from $x$ to
$y$ within time $t\geq0$ for a random walk on $\mathbb{Z}^{d}$
generated by matrix $A$. Set
$$G_{\lambda}(x,y):=\int\nolimits_{0}^{\infty}{e^{-\lambda
t}p(t;x,y)\,dt},\quad \lambda>0,\quad x,y\in\mathbb{Z}^{d}.$$ Note
that the Green's function
${G_{0}(x,y):=\lim\nolimits_{\lambda\to0+}{G_{\lambda}(x,y)}}$ is
well-defined and takes finite values for $d\geq3$ by virtue of the
transience of our random walk on $\mathbb{Z}^{d}$, $d\geq3$. One can
check (see \cite{VT_Siberia}) that ${h_{d}=(a G_{0}({\bf 0},{\bf
0}))^{-1}}$, $d\in\mathbb{N}$, where $a:=-a({\bf 0})$.

As shown in \cite{Y_Book}, Theorem 2.1.1 (see also
\cite{VT_Siberia}), for any fixed $x,y\in\mathbb{Z}^{d}$, one has
\begin{equation}\label{p(t;x,y)sim,p(t;0,0)-p(t;x,y)}
p(t;x,y)\sim\frac{\gamma_{d}}{t^{d/2}},\quad p(t;{\bf 0},{\bf
0})-p(t;x,y)\sim\frac{\tilde{\gamma}_{d}(y-x)}{t^{1+d/2}},\quad
t\to\infty,
\end{equation}
where
$\gamma_{d}:=\left((2\pi)^{d}\left|\det\phi''_{\theta\theta}({\bf
0})\right|\right)^{-1/2},$
$\phi(\theta):=\sum\nolimits_{z\in\mathbb{Z}^{d}}{a(z,{\bf
0})\cos(z,\theta)}$, $\theta\in[-\pi,\pi]^{d}$,
$$\phi''_{\theta\theta}({\bf
0})=\left(\left.\frac{\partial^{2}\phi(\theta)}{\partial
\theta_{i}\partial \theta_{j}}\right|_{\theta={\bf
0}}\right)_{i,j\in\{1,\ldots,d\}},\quad\tilde{\gamma}_{d}(z):=\frac{1}{2(2\pi)^{d}}\int\nolimits_{\mathbb{R}^{d}}{
(\upsilon,z)^{2}e^{(\phi''_{\theta\theta}({\bf
0})\upsilon,\upsilon)/2}\,d\upsilon},\quad z\in\mathbb{Z}^{d},$$ and
$(\cdot,\cdot)$ stands for the scalar product in $\mathbb{R}^{d}$.
In particular, it follows that the value
$$m_{d}:=1-(1-\alpha)a^{-1}+2(1-\alpha)a^{-1}G^{-2}_{0}({\bf 0},{\bf
0})\int\nolimits_{0}^{\infty}{t p(t;{\bf 0},{\bf 0})\,dt}$$ is
finite for $d\geq5$. Set also $q(s,t;x,y):=1-{\sf
E}_{x}{s^{\mu(t;y)}}$, $s\in[0,1]$, $t\geq0$,
$x,y\in\mathbb{Z}^{d}$. For $d=2$ we use the function
$$J(s;y):=\alpha\int\nolimits_{0}^{\infty}{(f(1-q(s,u;{\bf
0},y))-1+q(s,u;{\bf 0},y))\,d u},\quad s\in[0,1], \quad
y\in\mathbb{Z}^{d}.$$

The main results are contained in the following three theorems. For
the sake of completeness their statements include the case $y={\bf
0}$ studied earlier in \cite{B_TVP}--\cite{LMJ}, \cite{VT_4},
\cite{VT_Siberia}, \cite{VT_DM} and \cite{VTY}.

\begin{Th}\label{T:m(t;x,y)}
Let $x,y\in\mathbb{Z}^{d}$. The following relations are valid being
different for $y\neq{\bf 0}$ and $y={\bf 0}$, namely, as
$t\to\infty$,
\begin{eqnarray*}
m(t;x,y)\sim\frac{\gamma_{1}}{\sqrt{t}},\quad m(t;x,{\bf 0})\sim\frac{\gamma_{1}a}{(1-\alpha)\sqrt{t}},\quad&d=1,&\\
m(t;x,y)\sim\frac{\gamma_{2}}{t},\quad m(t;x,{\bf 0})\sim\frac{\gamma_{2}a}{(1-\alpha)t},\quad&d=2,&\\
m(t;x,y)\sim\frac{G_{0}(x,{\bf 0})G_{0}({\bf
0},y)}{2\pi\gamma_{3}\sqrt{t}},\quad m(t;x,{\bf 0})\sim\frac{a
G_{0}(x,{\bf 0})G_{0}({\bf 0},{\bf 0})}{2\pi\gamma_{3}(1-\alpha)\sqrt{t}},\quad&d=3,&\\
m(t;x,y)\sim\frac{G_{0}(x,{\bf 0})G_{0}({\bf
0},y)}{\gamma_{4}\,\ln{t}},\quad m(t;x,{\bf 0})\sim\frac{a
G_{0}(x,{\bf 0})G_{0}({\bf 0},{\bf 0})}
{\gamma_{4}(1-\alpha)\ln{t}},\quad&d=4,&\\
m(t;x,y)\to\frac{(1-\alpha)G_{0}(x,{\bf 0})G_{0}({\bf
0},y)}{a\,G^{2}_{0}({\bf 0},{\bf 0})\,m_{d}},\quad m(t;x,{\bf
0})\to\frac{G_{0}(x,{\bf 0})}{G_{0}({\bf 0},{\bf
0})m_{d}},\quad&d\geq5.&
\end{eqnarray*}
\end{Th}

\begin{Th}\label{T:q(t;x,y)sim}
For $x,y\in\mathbb{Z}^{d}$ and $t\to\infty$ the following formulae
hold true
\begin{eqnarray*}
q(t;x,y)\sim\frac{2(1-\alpha)}{\sigma^{2}\gamma_{1}\,a\sqrt{t}\ln{t}},\quad
&d=1,&\\
q(t;x,y)\sim\frac{\gamma_{2}}{t}\left(1-\frac{a}{1-\alpha}J(0;y)\right),\quad
y\neq{\bf 0},\quad &d=2,&\\
q(t;x,{\bf 0})\sim\frac{\gamma_{2}a}{(1-\alpha) t}(1-J(0;{\bf
0})),\quad &d=2,&\\
q(t;x,y)\sim\frac{4\pi\gamma_{3}(1-\alpha)G_{0}(x,{\bf
0})}{\sigma^{2}\,a\,G^{3}_{0}({\bf
0},{\bf 0})\sqrt{t}\ln{t}},\quad &d=3,&\\
q(t;x,y)\sim\frac{3\gamma_{4}(1-\alpha)G_{0}(x,{\bf
0})\ln{t}}{\sigma^{2}\,a\,G^{3}_{0}({\bf
0},{\bf 0})\,t},\quad &d=4,&\\
q(t;x,y)\sim\frac{2\,m_{d}G_{0}(x,{\bf 0})}{\sigma^{2}\,G_{0}({\bf
0},{\bf 0})t},\quad
&d\geq5,&\\
\end{eqnarray*}
where for $d=2$ and $s\in[0,1]$ the strict inequalities
$J(s;y)<(1-s)(1-\alpha)/a$, $y\neq{\bf 0}$, and $J(s;{\bf 0})<1-s$
are valid.
\end{Th}

\begin{Th}\label{T:limits}
Given $x,y\in\mathbb{Z}^{d}$, $\lambda\in[0,\infty)$ and
$s\in[0,1]$, one has, as $t\to\infty$,
\begin{eqnarray*}
\lim\limits_{t\to\infty}{\sf
E}_{x}\left(\left.\exp\left\{-\frac{\lambda\,\mu(t;y)}{{\sf
E}_{x}(\mu(t;y)|\mu(t;y)>0)}\right\}\right|\mu(t;y)>0\right)=\frac{1}{\lambda+1},\quad
d=1,\quad d=3\quad \mbox{or}&d\geq5,&\\
\lim\limits_{t\to\infty}{{\sf
E}_{x}\left(\left.s^{\mu(t;y)}\right|\mu(t;y)>0\right)}=\frac{(1-\alpha)s-a(J(0;y)-J(s;y))}{1-\alpha-a
J(0;y)},
\quad y\neq{\bf 0},\quad&d=2,&\\
\lim\limits_{t\to\infty}{{\sf E}_{x}\left(\left.s^{\mu(t;{\bf
0})}\right|\mu(t;{\bf 0})>0\right)}=\frac{s-(J(0;{\bf 0})-J(s;{\bf
0}))}{1-J(0;{\bf 0})},
\quad &d=2,&\\
\lim\limits_{t\to\infty}{{\sf
E}_{x}\left(\left.\exp\left\{-\frac{\lambda\,\mu(t;y)}{{\sf
E}_{x}(\mu(t;y)|\mu(t;y)>0)}\right\}\right|\mu(t;y)>0\right)}=\frac{1}{3}+\frac{2}{3}\cdot\frac{2}{2+3\lambda},\quad
&d=4.&
\end{eqnarray*}
\end{Th}

Observe that the normalizing factor ${\sf
E}_{x}{(\mu(t;y)|\mu(t;y)>0)}$ arising in Theorem \ref{T:limits} is
exactly $m(t;x,y)/q(t;x,y)$ and the asymptotic behavior of the
functions $m(t;x,y)$ and $q(t;x,y)$ is given by Theorems
\ref{T:m(t;x,y)} and \ref{T:q(t;x,y)sim}, respectively.

To establish Theorem \ref{T:m(t;x,y)} it is useful to invoke the
forward and backward Kolmogorov's differential equations (considered
in appropriate Banach spaces) for mean numbers of particles at
different points of the lattice and also the resulting integral
equations (see \cite{Y_Book}). As for Theorems \ref{T:q(t;x,y)sim}
and \ref{T:limits}, note that for proving results in \cite{B_TVP},
\cite{VT_4}, \cite{VT_Siberia}, \cite{VT_DM} and \cite{VTY}
concerning the number of particles at the origin the method of
introduction of an auxiliary Bellman-Harris branching process with
particles of \emph{two types} was efficient. However, for proving
Theorems \ref{T:q(t;x,y)sim} and \ref{T:limits} we have to involve a
Bellman-Harris branching process with particles of \emph{six types}.
To apply the latter method we attend to a new notion of the
\emph{hitting time with taboo} in the framework of a (non-branching)
random walk on $\mathbb{Z}^{d}$. More precisely, we use our recent
results (see \cite{B_Preprint}) on the asymptotic behavior of the
tail of the (improper) cumulative distribution function of this
time. Due to that one can employ the theorems by V.A.Vatutin for
Bellman-Harris branching processes with particles of several types
(see, e.g., \cite{Vatutin_1979}--\cite{V_1980}). Afterwards we have
to deal with sophisticated analytic estimates of the solutions of
the parametric integral equations (see, e.g., \cite{VT_4},
\cite{VT_Siberia}, \cite{VT_DM} and \cite{VTY}).

\section{Proof of Theorem \ref{T:m(t;x,y)}}

Let us recall some useful results employed within this section.
According to \cite{Gikhman_Skorokhod}, Ch.3, Sec.2, the transition
probabilities $p(t;x,y), t\geq0, x,y\in\mathbb{Z}^{d},$ of the
random walk generated by matrix $A$ satisfy the {\it backward}
Kolmogorov's equations
\begin{equation}\label{bKe_p(t;x,y)}
\frac{d\,p(t;x,y)}{d\,t}=\left(Ap(t;\cdot,y)\right)(x),\quad
p(0;x,y)=\delta_{y}(x).
\end{equation}
Here $(A
p(t;\cdot,y))(x)=\sum\nolimits_{z\in\mathbb{Z}^{d}}{a(x,z)p(t;z,y)}$
and $\delta_{y}(\cdot)$ is a column vector in the space
$l_{2}(\mathbb{Z}^{d})$ with zero components except for the
component 1 indexed by $y$. In a similar way, the {\it backward}
Kolmogorov's equations for $m(t;x,y)$, $t\geq0$,
$x,y\in\mathbb{Z}^{d}$, (see, e.g., Theorem 2.1 in \cite{Y_TVP})
take the form
\begin{equation}\label{bKe_m(t;x,y)}
\frac{d\,m(t;x,y)}{d\,
t}=\left(\overline{A}m(t;\cdot,y)\right)(x)+\overline{\beta}_{c}\left(\Delta_{\bf
0}m(t;\cdot,y)\right)(x),\quad m(0;x,y)=\delta_{y}(x),
\end{equation}
where
$\overline{A}=\left(\overline{a}(u,\upsilon)\right)_{u,\upsilon\in\mathbb{Z}^{d}}:=A+\left(a^{-1}(1-\alpha)-1\right)\Delta_{\bf
0}A$, $\Delta_{\bf 0}:=\delta_{\bf 0}\delta_{\bf 0}^{\textsc{T}}$
($\textsc{T}$ stands for transposition) and
$\overline{\beta}_{c}:=(1-\alpha)a^{-1}G^{-1}_{0}({\bf 0},{\bf 0})$.
Here we follow the notation of \cite{Y_TVP}.

\begin{Lm}\label{L:monotonicity_m(t;y,y)}
For each $y\in\mathbb{Z}^{d}$, the function $m(t;y,y)$ is
non-increasing in $t$.
\end{Lm}
{\sc Proof. }The monotonicity of $m(\cdot;y,y)$ for SBRW on
$\mathbb{Z}^{d}$ was established in Lemma 3.3.5 of \cite{Y_Book}.
The key step of its proof was to use self-adjointness of the
operator $H:=A+\beta_{c}\Delta_{\bf 0}$ where
$\beta_{c}:=G^{-1}_{0}({\bf 0},{\bf 0})$. For CBRW the analog of $H$
is the non self-adjoint operator
$\overline{H}:=\overline{A}+\overline{\beta}_{c}\Delta_{\bf 0}$.
However, Lemma 3.1 in \cite{Y_TVP} permits to pass to (self-adjoint)
{\it symmetrization} of $\overline{H}$ and then apply Lemma 3.3.5 in
\cite{Y_Book}. Further argument is similar to the proof of Theorem
in \cite{Y_Vestnik}. $\square$

Equation \eqref{bKe_m(t;x,y)} was obtained by differentiating at
$s=1$ the following {\it backward} Kolmogorov's equation for the
generating function $F(s,t;x,y):={\sf E}_{x}{s^{\mu(t;y)}}$,
$s\in[0,1]$, $t\geq0$, $x,y\in\mathbb{Z}^{d}$, (see~\cite{Y_TVP})
\begin{equation}\label{bKe_F(s,t;x,y)}
\frac{\partial F(s,t;x,y)}{\partial
t}=\left(\overline{A}F(s,t;\cdot,y)\right)(x)+\left(\Delta_{\bf
0}\overline{f}(F(s,t;\cdot,y))\right)(x),\quad
F(s,0;x,y)=s^{\delta_{y}(x)}.
\end{equation}
Here $\overline{f}(s):=\alpha(f(s)-s)$, $s\in[0,1]$, is an
infinitesimal generating function of the number of offsprings of a
parent particle. We will employ \eqref{bKe_F(s,t;x,y)} in Section 5.

In Lemma \ref{L:fKe_F(s,t;x,y)} we derive a counterpart of the
forward Kolmogorov's equation for the fuction $F(s,t;x,y)$,
$s\in[0,1]$, $t\geq0$, $x,y\in\mathbb{Z}^{d}$, and, as a
consequence, the forward Kolmogorov's equation for $m(t;x,y)$.
Recall that $\overline{A}^{\ast}$ denotes an adjoint operator for
$\overline{A}$ and
$\left(\overline{A}^{\ast}m(t;x,\cdot)\right)(y)=\sum\nolimits_{z\in\mathbb{Z}^{d}}{m(t;x,z)\overline{a}(z,y)}$.

\begin{Lm}\label{L:fKe_F(s,t;x,y)}
For $s\in[0,1]$, $t\geq0$, $x,y\in\mathbb{Z}^{d}$, the following
relation holds true
\begin{eqnarray}
\frac{\partial F(s,t;x,y)}{\partial
t}&=&(s-1)\sum_{z\in\mathbb{Z}^{d},\,z\neq y}{\overline{a}(z,y){\sf
E}_{x}{s^{\mu(t;y)}\mu(t;z)}}+(s-1)\overline{a}(y,y){\sf
E}_{x}{s^{\mu(t;y)-1}\mu(t;y)}\nonumber\\
&+&\delta_{\bf 0}(y)\overline{f}(s){\sf
E}_{x}{s^{\mu(t;y)-1}\mu(t;y)},\quad
F(s,0;x,y)=s^{\delta_{x}(y)}.\label{fKe_F(s,t;x,y)}
\end{eqnarray}
Moreover, one has
\begin{equation}\label{fKe_m(t;x,y)}
\frac{d\,m(t;x,y)}{d\,t}=\left(\overline{A}^{\ast}m(t;x,\cdot)\right)(y)+\overline{\beta}_{c}\left(\Delta_{\bf
0}m(t;x,\cdot)\right)(y),\quad m(0;x,y)=\delta_{x}(y).
\end{equation}
\end{Lm}
{\sc Proof. }As usual in derivation of forward Kolmogorov's
equations, we consider all possible evolutions of the particles
population within the time interval $[t,t+h)$ and let $h\to0+$. To
justify arising passages to the limit we involve the Lebesgue
theorem on dominated convergence and useful estimates for transition
probabilities (see proof of Lemma 3 in \cite{Gikhman_Skorokhod},
Ch.3, Sec.2). We also benefit from finiteness of the mean total
number of particles $M(t;x):={\sf
E}_{x}\left(\sum\nolimits_{z\in\mathbb{Z}^{d}}{\mu(t;z)}\right)$ for
each $x\in\mathbb{Z}^{d}$ and $t\geq0$. The latter observation is
true since the last function belonging to
$l_{\infty}(\mathbb{Z}^{d})$ is a solution of the linear
differential equation in \eqref{bKe_m(t;x,y)} with the initial
condition ${M(0;x)=1}$ for all $x$ (instead of $\delta_{y}(x)$ in
\eqref{bKe_m(t;x,y)}), see \cite{Y_TVP}.

Equation \eqref{fKe_m(t;x,y)} is an immediate consequence of
\eqref{fKe_F(s,t;x,y)} due to formula
$m(t;x,y)=\partial_{s}{F(s,t;x,y)}|_{s=1}$. We also take into
account that $\overline{f}\,'(1)=\overline{\beta}_{c}$ in view of
\eqref{assumption_3}. $\square$

Consider equations \eqref{bKe_m(t;x,y)} and \eqref{fKe_m(t;x,y)} as
inhomogeneous ones for differential equation \eqref{bKe_p(t;x,y)} in
Banach space $l_{\infty}(\mathbb{Z}^{d})$. Applying the variation of
constant formula (see \cite{D_Krein}, Ch.2, Sec.1) we infer that
\begin{eqnarray}
m(t;x,y)=p(t;x,y)&+&\left(1-\frac{a}{1-\alpha}\right)\int\nolimits_{0}^{t}{p(t-u;x,{\bf
0})m'(u;{\bf
0},y)\,du}\nonumber\\
&+&\frac{a\,\overline{\beta}_{c}}{1-\alpha}\int\nolimits_{0}^{t}{p(t-u;x,{\bf
0})m(u;{\bf 0},y)\,du},\label{int_bKe_m(t;x,y)}\\
m(t;x,y)=p(t;x,y)&+&\left(\frac{1-\alpha}{a}-1\right)\int\nolimits_{0}^{t}{m(t-u;x,{\bf
0})p'(u;{\bf
0},y)\,du}\nonumber\\
&+&\overline{\beta}_{c}\int\nolimits_{0}^{t}{m(t-u;x,{\bf
0})p(u;{\bf 0},y)\,du}.\label{int_fKe_m(t;x,y)}
\end{eqnarray}
An analogous result for SBRW on $\mathbb{Z}^{d}$ can be found in
\cite{Y_Book}, Theorem 1.4.1. Now we can give\\
{\sc Proof of Theorem \ref{T:m(t;x,y)}. }To find the asymptotic
behavior of $m(t;x,y)$, $t\to\infty$, $x,y\in\mathbb{Z}^{d}$,
$y\neq{\bf 0}$, we estimate each of the summands in the right-hand
sides of \eqref{int_bKe_m(t;x,y)} and \eqref{int_fKe_m(t;x,y)} when
$x\neq{\bf 0}$ and $x={\bf 0}$, respectively, as $t\to\infty$.
Namely, we will show that, for $d=1$ and $d=2$, the main
contribution to the asymptotic behavior of the right-hand side of
\eqref{int_bKe_m(t;x,y)}, as well as of \eqref{int_fKe_m(t;x,y)}, is
due to the first summand. However, for $d\geq3$, the asymptotic
behavior of the right-hand sides of \eqref{int_bKe_m(t;x,y)} and
\eqref{int_fKe_m(t;x,y)} is determined only by the third summands.
It is worth mentioning that, for $d=1$ and $d=2$, the third summands
in \eqref{int_bKe_m(t;x,y)} and \eqref{int_fKe_m(t;x,y)} vanish in
view of equality $\overline{\beta}_{c}=0$.

Let $x={\bf 0}$. The asymptotic behavior of the first summand in the
right-hand side of \eqref{int_fKe_m(t;x,y)} is given by
\eqref{p(t;x,y)sim,p(t;0,0)-p(t;x,y)}. The estimate of the second
summand could be obtained on account of Lemma~6 in \cite{VT_Siberia}
and, in particular, relation (20). However, to avoid verifying the
bounded variation of the functions $p(t;x,y)$ and $m(t;x,y)$ in
variable $t$ we choose another approach consisting in direct
estimation of the second summand. Recall that representation
(2.1.15) in \cite{Y_Book} entails the inequalities $p'(t;{\bf
0},{\bf 0})\leq0$, $p'(t;{\bf 0},{\bf 0})\leq p'(t;{\bf 0},y)$,
$p''(t;{\bf 0},{\bf 0})\geq0$ and ${p''(t;{\bf 0},{\bf
0})-p''(t;{\bf 0},y)\geq0}$, $t\geq0$. Then by virtue of
\eqref{p(t;x,y)sim,p(t;0,0)-p(t;x,y)} as well as the classical
results on differentiating the asymptotic formulae (see, e.g.,
\cite{Brein}, Ch.7, Sec.3), for $d\in\mathbb{N}$, one has
$$p'(t;{\bf 0},{\bf 0})\sim-\frac{d\,\gamma_{d}}{2\,t^{d/2+1}},\quad
p'(t;{\bf 0},{\bf 0})-p'(t;{\bf
0},y)\sim-\frac{(d+2)\tilde{\gamma}_{d}(y)}{2t^{d/2+2}},\quad
t\to\infty.$$ Whence taking into account Lemma 5.1.2 in
\cite{Y_Book} ("lemma on convolutions") and the already proved
assertion of Theorem \ref{T:m(t;x,y)} for $x={\bf 0}$ and $y={\bf
0}$ we deduce that, as $t\to\infty$,
\begin{eqnarray}
\!\!\!\!& &\int\nolimits_{0}^{t}{m(t-u;{\bf 0},{\bf 0})p'(u;{\bf
0},y)\,du} =\int\nolimits_{0}^{t}{m(t-u;{\bf 0},{\bf
0})\left(p'(u;{\bf 0},y)-p'(u;{\bf
0},{\bf 0})\right)\,du}\nonumber\\
\!\!\!\!\!\!&+&\!\!\!\!\int\nolimits_{0}^{t}{m(t-u;{\bf 0},{\bf
0})p'(u;{\bf 0},{\bf 0})\,du}=m(t;{\bf 0},{\bf 0})-m(t;{\bf 0},{\bf
0})+o(m(t;{\bf 0},{\bf 0}))=o(m(t;{\bf 0},{\bf
0})).\quad\quad\quad\label{second_summand}
\end{eqnarray}
Combining relations \eqref{p(t;x,y)sim,p(t;0,0)-p(t;x,y)},
\eqref{int_fKe_m(t;x,y)} and \eqref{second_summand} we establish
Theorem \ref{T:m(t;x,y)} for $d=1$ and $d=2$ when $x={\bf 0}$. The
statement of Theorem \ref{T:m(t;x,y)} for $d\geq3$ and $x={\bf 0}$
follows from formulae \eqref{p(t;x,y)sim,p(t;0,0)-p(t;x,y)},
\eqref{int_fKe_m(t;x,y)} and \eqref{second_summand} by Lemma 5.1.2
in \cite{Y_Book} and in view of Theorem \ref{T:m(t;x,y)} for the
known case $x=y={\bf 0}$.

Let $x\neq{\bf 0}$. Similarly to the case $x={\bf 0}$ we see that
\begin{equation}\label{second_summand_x}
\int\nolimits_{0}^{t}{p(t-u;x,{\bf 0})m'(u;{\bf
0},y)\,du}=\int\nolimits_{0}^{t}{m(t-u;{\bf 0},y)p'(u;x,{\bf
0})\,du}=o(m(t;{\bf 0},y)),\quad t\to\infty.
\end{equation}
Thus, the combination of \eqref{p(t;x,y)sim,p(t;0,0)-p(t;x,y)},
\eqref{int_bKe_m(t;x,y)} and \eqref{second_summand_x} proves Theorem
\ref{T:m(t;x,y)} for $d=1$ or $d=2$ and $x\neq{\bf 0}$. For $d\geq3$
and $x\neq{\bf 0}$ we estimate the third summand in
\eqref{int_bKe_m(t;x,y)} with the help of Lemma 5.1.2 in
\cite{Y_Book}, relation \eqref{p(t;x,y)sim,p(t;0,0)-p(t;x,y)} and
the assertion of Theorem \ref{T:m(t;x,y)} for $d\geq3$ and $x={\bf
0}$ established above. $\square$

\section{Auxiliary Bellman-Harris branching process}

Let us briefly describe a Bellman-Harris branching process with
particles of six types. It is initiated by a single particle of type
$i=1,\ldots,6$. The parent particle has a random life-length with a
cumulative distribution function (c.d.f.) $G_{i}(t)$, $t\geq0$. When
dying the particle produces offsprings according to a generating
function $f_{i}(\vec{s}\,),$
$\vec{s}=(s_{1},\ldots,s_{6})\in[0,1]^{6}$. The new particles of
type $j=1,\ldots,6$ evolve independently with the life-length
distribution $G_{j}(t)$ and an offspring generating function
$f_{j}(\vec{s}\,)$. Let $M:=\left(\partial_{s_{j}}
f_{i}|_{\vec{s}=(1,\ldots,1)}\right)_{i,j=1,\ldots,6}$ be the mean
matrix of the process. The Bellman-Harris branching process is
called {\it critical indecomposable} if the Perron root of $M$ (i.e.
eigenvalue having the maximal modulus) equals 1 and for some integer
$n$ all elements of $M^{n}$ are positive (see, e.g., \cite{Sev},
Ch.4, Sec.6 and 7). Denote the number of particles of type $j$
existing at time $t$ by $Z_{j}(t)$, $t\geq0$, $j=1,\ldots,6$. Set
$F_{i}(t;\vec{s}\,)={\sf
E}_{i}\left(\prod\nolimits_{j=1}^{6}{s^{Z_{j}(t)}_{j}}\right)$,
$i=1,\ldots,6$, $t\geq0$, $\vec{s}\in[0,1]^{6}$, where the index $i$
means that the parent particle is of type $i$. In other words,
$F_{i}(t;\vec{s}\,)$ is a generating function of the numbers of
particles of all types existing at time $t$ given that the process
is initiated by a single particle of type $i$.

Before demonstrating how an auxiliary Bellman-Harris process can be
constructed in the framework of CBRW on $\mathbb{Z}^{d}$ we have to
introduce some notation. Recall that in \cite{B_Preprint} a new
notion of a hitting time with taboo was proposed for a
(non-branching) random walk on $\mathbb{Z}^{d}$ generated by matrix
$A$. Namely, let $\tau^{-}_{y,z}$, $y,z\in\mathbb{Z}^{d}$, $y\neq
z$, be the time spent by the particle (performing the random walk)
after leaving the starting point until the first hitting $y$ if
particle's trajectory does not pass $z$. Otherwise (if particle's
trajectory passes point $z$ before the first hitting $y$),
$\tau^{-}_{y,z}=\infty$. Denote by $H^{-}_{x,y,z}(t)$, $t\geq0$, the
improper c.d.f. of $\tau^{-}_{y,z}$ given that the starting point of
the random walk is $x\in\mathbb{Z}^{d}$.

Return to CBRW on $\mathbb{Z}^{d}$. In this section we assume that
CBRW may start at the origin or at a fixed point $y\neq{\bf 0}$. We
divide the particles population existing at time $t\geq0$ into seven
groups. The particles located at time $t$ at the origin
(respectively, at $y$) form the first (respectively, second) group
having cardinality $\mu(t;{\bf 0})$ (respectively, $\mu(t;y)$). Next
consider at time $t$ a family of particles labeled by a collection
$(u,v,w)$ of lattice points, its cardinality being $\mu_{u,v,w}(t)$.
It consists of the particles which have left $u$ at least once
within time interval $[0,t]$, upon the last leaving $u$ have yet
reached neither $v$ nor $w$ but eventually will hit $v$ before
possible hitting $w$. Our third group corresponds to $(u,v,w)=({\bf
0},y,{\bf 0})$, the fourth to $(y,{\bf 0},y)$, the fifth to $({\bf
0},{\bf 0},y)$ and the sixth to $(y,y,{\bf 0})$. The seventh group
comprises the rest of particles not included into the above six
groups. Note that the last group consists of the particles having
infinite life-length since after time $t$ they will not hit the
origin any more. So, after time $t$ these particles will not produce
any offsprings and have no influence on the numbers of particles in
other six groups.

Now we can introduce an auxiliary Bellman-Harris process and use it
for the study of CBRW on $\mathbb{Z}^{d}$. Consider a
six-dimensional Bellman-Harris process having the following c.d.f.
$G_{i}$ and generating function $f_{i}$, $i=1,\ldots,6$,
\begin{eqnarray*}
G_{1}(t)\!&=&\!1-e^{-t},\quad f_{1}(\vec{s}\,)=\alpha
f(s_{1})+(1-\alpha)H^{-}_{{\bf 0},y,{\bf
0}}(0)s_{2}+(1-\alpha)(H^{-}_{{\bf 0},y,{\bf 0}}(\infty)-H^{-}_{{\bf
0},y,{\bf 0}}(0))s_{3}\\
& &\quad\quad\quad\quad\quad\quad\quad+(1-\alpha)H^{-}_{{\bf 0},{\bf
0},y}(\infty)s_{5}+(1-\alpha)(1-H^{-}_{{\bf 0},y,{\bf
0}}(\infty)-H^{-}_{{\bf 0},{\bf 0},y}(\infty)),\\
G_{2}(t)\!&=&\!1-e^{-a t},\quad f_{2}(\vec{s}\,)=H^{-}_{y,{\bf
0},y}(0)s_{1}+(H^{-}_{y,{\bf 0},y}(\infty)-H^{-}_{y,{\bf
0},y}(0))s_{4}\\
& &\quad\quad\quad\quad\quad\quad\quad\;+H^{-}_{y,y,{\bf
0}}(\infty)s_{6}+(1-H^{-}_{y,{\bf
0},y}(\infty)-H^{-}_{y,y,{\bf 0}}(\infty)),\\
G_{3}(t)\!&=&\!\frac{H^{-}_{{\bf 0},y,{\bf 0}}(t)-H^{-}_{{\bf
0},y,{\bf 0}}(0)}{H^{-}_{{\bf 0},y,{\bf 0}}(\infty)-H^{-}_{{\bf
0},y,{\bf 0}}(0)},\quad f_{3}(\vec{s}\,)=s_{2},\quad
G_{4}(t)=\frac{H^{-}_{y,{\bf 0},y}(t)-H^{-}_{y,{\bf
0},y}(0)}{H^{-}_{y,{\bf
0},y}(\infty)-H^{-}_{y,{\bf 0},y}(0)},\quad f_{4}(\vec{s}\,)=s_{1},\\
G_{5}(t)\!&=&\!\frac{H^{-}_{{\bf 0},{\bf 0},y}(t)}{H^{-}_{{\bf
0},{\bf 0},y}(\infty)},\quad f_{5}(\vec{s}\,)=s_{1},\quad
G_{6}(t)=\frac{H^{-}_{y,y,{\bf 0}}(t)}{H^{-}_{y,y,{\bf
0}}(\infty)},\quad f_{6}(\vec{s}\,)=s_{2},
\end{eqnarray*}
where
$H^{-}_{x,y,z}(\infty):=\lim\nolimits_{t\to\infty}{H^{-}_{x,y,z}(t)}$.
The symmetry and homogeneity of the random walk generated by matrix
$A$ imply identities $H^{-}_{{\bf 0},y,{\bf 0}}\equiv H^{-}_{y,{\bf
0},y}$ and $H^{-}_{{\bf 0},{\bf 0},y}\equiv H^{-}_{y,y,{\bf 0}}$,
whence $G_{3}\equiv G_{4}$ and $G_{5}\equiv G_{6}$. It is not
difficult to see that for the branching process constructed in this
way one has $(\mu(t;{\bf 0}),\mu(t;y),\mu_{{\bf 0},y,{\bf
0}}(t),\mu_{y,{\bf 0},y}(t),\mu_{{\bf 0},{\bf 0},y}(t),\mu_{y,y,{\bf
0}}(t))\stackrel{Law}=(Z_{1}(t),\ldots,Z_{6}(t))$, $t\geq0$.

Observe that the introduced Bellman-Harris branching process with
particles of six types is critical indecomposable. Indeed, it is an
easy computation task to check that all entries of $M^{6}$ are
positive. Furthermore, if $H^{-}_{{\bf 0},y,{\bf 0}}(0)\neq0$ (that
is $a({\bf 0},y)>0$) then already all entries of $M^{4}$ are
positive. Hence, the constructed process is indecomposable. To
verify its criticality note that in view of Theorem 3 in
\cite{B_Preprint} one can rewrite the first relation in
\eqref{assumption_3} as follows
\begin{equation}\label{criticality_condition_new}
\alpha f'(1)=1-(1-\alpha)\left(H^{-}_{{\bf 0},{\bf
0},y}(\infty)+\frac{(H^{-}_{{\bf 0},y,{\bf
0}}(\infty))^{2}}{1-H^{-}_{{\bf 0},{\bf 0},y}(\infty)}\right).
\end{equation}
Then by inspecting the explicit expression for the characteristic
polynomial of the mean matrix $M$ we deduce that it has the form
$$\det{(M-\kappa I)}=\kappa^{2}(\kappa-1)R(\kappa)$$
where $I$ is a unit matrix, $\kappa\in\mathbb{C}$ and
\begin{eqnarray*}
R(\kappa):=\kappa^{3}&+&\kappa^{2}(1-\alpha
f'(1))+\kappa\left(1-\alpha f'(1)-(2-\alpha)H^{-}_{{\bf 0},{\bf
0},y}(\infty)-(1-\alpha)(H^{-}_{{\bf 0},y,{\bf
0}}(0))^{2}\right)\\
&+&(1-\alpha)(H^{-}_{{\bf 0},y,{\bf 0}}(\infty)-H^{-}_{{\bf
0},y,{\bf 0}}(0))^{2}-(1-\alpha)(H^{-}_{{\bf 0},{\bf
0},y}(\infty))^{2}.
\end{eqnarray*}
The polynomial $R(\kappa)$ has no real roots greater than 1 because
$R(1)>0$ and $R\,'(\kappa)>0$ for $\kappa\geq1$. In fact, due to
identity \eqref{criticality_condition_new} we obtain the
representation with strictly positive summands
\begin{eqnarray*}
R(1)&=&(1-H^{-}_{{\bf 0},{\bf 0},y}(\infty))((1-\alpha)H^{-}_{{\bf
0},{\bf 0},y}(\infty)+1)+(1-\alpha)H^{-}_{{\bf 0},y,{\bf
0}}(\infty)(H^{-}_{{\bf 0},y,{\bf 0}}(\infty)-H^{-}_{{\bf 0},y,{\bf
0}}(0))\\
&+&\frac{(1-\alpha)H^{-}_{{\bf 0},y,{\bf 0}}(\infty)(H^{-}_{{\bf
0},y,{\bf 0}}(\infty)-H^{-}_{{\bf 0},y,{\bf 0}}(0)(1-H^{-}_{{\bf
0},{\bf 0},y}(\infty)))}{1-H^{-}_{{\bf 0},{\bf
0},y}(\infty)}+\frac{(1-\alpha)(H^{-}_{{\bf 0},y,{\bf
0}}(\infty))^{2}}{1-H^{-}_{{\bf 0},{\bf 0},y}(\infty)}.
\end{eqnarray*}
Moreover, if $\kappa\geq1$ then
\begin{eqnarray*}
R\,'(\kappa)&=&3\kappa^{2}+2\kappa(1-\alpha f'(1))+1-\alpha
f'(1)-(2-\alpha)H^{-}_{{\bf 0},{\bf
0},y}(\infty)-(1-\alpha)(H^{-}_{{\bf 0},y,{\bf
0}}(0))^{2}\\
&>&3-2H^{-}_{{\bf 0},{\bf 0},y}(\infty)-H^{-}_{{\bf 0},y,{\bf
0}}(0)>0.
\end{eqnarray*}
Thus, the greatest positive real root of the characteristic
polynomial of $M$ is 1. Hence, by the Frobenius theorem (see, e.g.,
Theorem 2 in \cite{Sev}, Ch.4, Sec.5) 1 is the Perron root of $M$.
So, the auxiliary Bellman-Harris process is critical.

Denote by $\vec{v}=(v_{1},\ldots,v_{6})$ and
$\vec{u}=(u_{1},\ldots,u_{6})$ the left and right positive
eigenvectors corresponding to the Perron root of $M$ such that
$(\vec{u},\vec{1}\,)=1$ and $(\vec{v},\vec{u}\,)=1$ where
$\vec{1}=(1,\ldots,1)\in\mathbb{R}^{6}$. Taking into account
\eqref{criticality_condition_new} we rewrite the components of
$\vec{u}$ and $\vec{\upsilon}$ in the convenient form
\begin{eqnarray}
u_{1}=u_{4}=u_{5}=\frac{1-H^{-}_{{\bf 0},{\bf 0},y}(\infty)}{U},& &
u_{2}=u_{3}=u_{6}=\frac{H^{-}_{{\bf
0},y,{\bf 0}}(\infty)}{U},\label{u=}\\
v_{1}=\frac{U}{V},\quad v_{2}=\frac{U(1-\alpha)H^{-}_{{\bf 0},y,{\bf
0}}(\infty)}{V(1-H^{-}_{{\bf 0},{\bf 0},y}(\infty))},&
&v_{3}=\frac{U(1-\alpha)(H^{-}_{{\bf 0},y,{\bf
0}}(\infty)-H^{-}_{{\bf 0},y,{\bf
0}}(0))}{V},\label{v=1}\\
v_{4}=v_{2}(H^{-}_{{\bf 0},y,{\bf 0}}(\infty)-H^{-}_{{\bf 0},y,{\bf
0}}(0)),& & v_{5}=\frac{U(1-\alpha)H^{-}_{{\bf 0},{\bf
0},y}(\infty)}{V},\quad v_{6}=v_{2}H^{-}_{{\bf 0},{\bf
0},y}(\infty)\quad\label{v=2}
\end{eqnarray}
where the auxiliary variables $U$ and $V$ are defined by way of
$$U:=3(1-H^{-}_{{\bf 0},{\bf 0},y}(\infty)+H^{-}_{{\bf 0},y,{\bf 0}}(\infty)),$$
$$V:=3-2\alpha f'(1)-(2-\alpha)H^{-}_{{\bf 0},{\bf
0},y}(\infty)+(1-\alpha)((H^{-}_{{\bf 0},y,{\bf
0}}(\infty)-H^{-}_{{\bf 0},y,{\bf 0}}(0))^{2}-(H^{-}_{{\bf 0},y,{\bf
0}}(0))^{2}-(H^{-}_{{\bf 0},{\bf 0},y}(\infty))^{2}).$$

Using decomposition $f(1-x)=1-f'(1)x+f''(1)x^{2}/2+o(x^{2}),$
$x\to0+$, along with formulae
\eqref{criticality_condition_new}--\eqref{v=2} and the definition of
$\vec{f}(\vec{s}\,)=(f_{1}(\vec{s}\,),\ldots,f_{6}(\vec{s}\,))$, it
is not difficult to verify by standard calculations that
\begin{equation}\label{decomposition}
x-\left(\vec{v},\vec{1}-\vec{f}\left(\vec{1}-\vec{u}x\right)\right)\sim
B x^{2},\quad x\to0+,\quad\mbox{where}\quad
B:=\frac{\sigma^{2}\left(1-H^{-}_{{\bf 0},{\bf
0},y}(\infty)\right)^{2}}{2 UV}.
\end{equation}

In the next two lemmas we apply theorems proved in papers
\cite{Vatutin_1979}--\cite{V_1980} to the constructed
six-dimensional Bellman-Harris branching process and then
reformulate the obtained results for CBRW on $\mathbb{Z}^{d}$ when
$d\geq5$. Common to these theorems are the conditions of criticality
and indecomposability of the Bellman-Harris process which were
established above. Another common condition on the behavior of the
function $x-(\vec{v},\vec{1}-\vec{f}(\vec{1}-\vec{u}x))$ is
fulfilled due to \eqref{decomposition}. However, various Vatutin's
theorems involve different assumptions on the order of asymptotic
decrease of the tails of $G_{k}(\cdot)$, $k=1,\ldots,6$. It is worth
to mention that such asymptotic behavior was established in
\cite{B_Preprint}, Theorem 3. Namely, our result for $d\leq5$
corresponds to condition of Theorem 1 in \cite{Vatutin_1979} whereas
the cases $d=6$ and $d\geq7$ meet the respective conditions of
Theorem~3 in \cite{V_1980} and Theorem 2 in \cite{Vatutin_1978}.

\begin{Lm}\label{L:d=5}
Given $y\in\mathbb{Z}^{5}$, $y\neq{\bf 0}$, for CBRW on
$\mathbb{Z}^{5}$ one has
$$q(t;{\bf 0},y)=o\left(t^{-3/4}\right),\quad q(t;y,y)=o\left(t^{-3/4}\right),\quad t\to\infty.$$
\end{Lm}
{\sc Proof. } To apply Theorem 1 in \cite{Vatutin_1979} to the
six-dimensional Bellman-Harris process constructed above for CBRW on
$\mathbb{Z}^{5}$ we verify the conditions of that theorem. According
to the definition of
$\vec{G}(\cdot)=(G_{1}(\cdot),\ldots,G_{6}(\cdot))$ and by Theorem 3
in \cite{B_Preprint} for $d=5$, the variable $\beta$ in condition 2)
of Theorem 1 in \cite{Vatutin_1979} is equal to $3/2$ whereas the
function $L_{1}(t)$ in the same condition tends to a constant, as
$t\to\infty$. The validity of condition 3) of Theorem 1 in
\cite{Vatutin_1979} is implied by Theorem~1 in \cite{Vatutin_1977}
(for our process the function $L_{1}(t)$ in this theorem tends to
$1/B$, as $t\to\infty$, in view of \eqref{decomposition}) combined
with the definition of $\vec{G}(\cdot)$ and Theorem 3 in
\cite{B_Preprint} for $d=5$. Thus, we may employ Theorem~1 in
\cite{Vatutin_1979}. Taking into account Theorem 3 in
\cite{B_Preprint} for $d=5$ and formulae \eqref{u=}--\eqref{v=2} we
deduce from Theorem 1 in \cite{Vatutin_1979} that
$\lim\nolimits_{t\to\infty}{\sf
E}_{i}(s_{2}^{Z_{2}(t)}|\vec{Z}(t)\neq\vec{0}\,)=1$ for each
$s_{2}\in[0,1]$ and $i=1,2$ (as usual,
$\vec{Z}(t)=(Z_{1}(t),\ldots,Z_{6}(t))$ and
$\vec{0}=(0,\ldots,0)\in\mathbb{R}^{6}$). Setting $s_{2}=0$ in the
last relation one has ${\sf P}_{i}(Z_{2}(t)>0)=o({\sf
P}_{i}(\vec{Z}(t)\neq\vec{0}\,))$, as $t\to\infty$. Moreover,
examining the proof of Theorem 1 in \cite{Vatutin_1979} we can show
that for our Bellman-Harris process the slowly varying function
$L^{\ast}(x)$ in the assertion of that theorem turns equivalent to
$1/\sqrt{B}$, as $x\to0+$. Consequently, the indicated in
\cite{Vatutin_1979} formula (0.4) can be sharpened in our case,
namely, the function ${\sf P}_{i}(\vec{Z}(t)\neq\vec{0}\,)$ has an
order of decreasing $t^{-3/4}$, as $t\to\infty$. Whence by the
connection between CBRW on $\mathbb{Z}^{5}$ and the auxiliary
Bellman-Harris process we complete the proof. $\square$

\begin{Lm}\label{L:d>=6}
In the framework of CBRW on $\mathbb{Z}^{d}$ with $d\geq6$ the
following relations hold true for $y\in\mathbb{Z}^{d}$, $y\neq{\bf
0}$,
$$q(t;{\bf 0},y)\sim\frac{2\,m_{d}}{\sigma^{2}\,t},\quad
q(t;y,y)\sim\frac{2\,m_{d}\,G_{0}({\bf 0},y)}{\sigma^{2}\,G_{0}({\bf 0},{\bf 0})\,t},\quad t\to\infty.$$
\end{Lm}
{\sc Proof. }Let us apply Theorem 3 in \cite{V_1980} to our
Bellman-Harris branching process when $d=6$. To this end we verify
whether all the conditions of Theorem 3 in \cite{V_1980} are
satisfied. In view of \eqref{decomposition} relation (6) in
\cite{V_1980} is valid for our process and the function $L_{1}(n)$
in (6) tends to $1/B$, as $n\to\infty$. Equality (7) in
\cite{V_1980} is also satisfied due to (6) in \cite{V_1980} in view
of the definition of $\vec{G}(\cdot)$ and Theorem~3 in
\cite{B_Preprint} for $d=6$. Now we may apply Theorem 3 in
\cite{V_1980}. In particular, it follows that for each $i=1,2$ the
expressions $\lim\nolimits_{t\to\infty}{{\sf
P}_{i}(Z_{1}(t)=0|\vec{Z}(t)\neq\vec{0}\,)}$ and
$\lim\nolimits_{t\to\infty}{{\sf
P}_{i}(Z_{2}(t)=0|\vec{Z}(t)\neq\vec{0}\,)}$ coincide, are positive
and strictly less than 1. Consequently, ${\sf
P}_{i}(Z_{1}(t)>0)\sim{\sf P}_{i}(Z_{2}(t)>0)$, as $t\to\infty$. The
asymptotic behavior of $q(t;{\bf 0},{\bf 0})={\sf
P}_{1}(Z_{1}(t)>0)$ and $q(t;y,{\bf 0})={\sf P}_{2}(Z_{1}(t)>0)$ can
be found in \cite{LMJ}, Lemmas 2 and 4. Thus, Lemma \ref{L:d>=6} is
proved for $d=6$.

For $d\geq7$ we will employ Theorem 2 in \cite{Vatutin_1978}.
Condition (6) of that theorem is valid due to Theorem 1 in
\cite{Vatutin_1977} (for our process, the function $L_{1}(t)$ in
this theorem tends to $1/B$, as $t\to\infty$) by virtue of the
definition of $\vec{G}(\cdot)$ and Theorem 3 in \cite{B_Preprint}
for $d\geq7$. The definition of $G_{k}(\cdot)$ and Theorem 3 in
\cite{B_Preprint} for $d\geq7$ also imply that
$\int\nolimits_{0}^{\infty}{t\,d G_{k}}<\infty$ for each
$k=1,\ldots,6$. So, all the conditions of Theorem 2 in
\cite{Vatutin_1978} are satisfied and it follows that
${\lim\nolimits_{t\to\infty}{{\sf
P}_{i}(Z_{k}(t)=0|\vec{Z}(t)\neq\vec{0}\,)}=0}$ for each
$k=1,\ldots,6$ and $i=1,2$. Hence, ${\sf P}_{i}(Z_{1}(t)>0)\sim{\sf
P}_{i}(Z_{2}(t)>0),$ $t\to\infty$. Notably, the asymptotic behavior
of $q(t;{\bf 0},{\bf 0})={\sf P}_{1}(Z_{1}(t)>0)$ and $q(t;y,{\bf
0})={\sf P}_{2}(Z_{1}(t)>0)$ can be found in \cite{LMJ}, Lemmas 2
and 4. Lemma \ref{L:d>=6} is proved for $d\geq7$. $\square$

Concluding this section we derive an integral equation in function
$q(\cdot;{\bf 0},y)$, ${y\neq{\bf 0}}$, which is a counterpart of
equation (2.6) in \cite{VT_4} for $q(\cdot;{\bf 0},{\bf 0})$. Our
integral equation will be essentially used for proving Theorem
\ref{T:q(t;x,y)sim} when $d=4$. Before formulating the corresponding
statement we have to introduce some more notation. Let $\tau_{z}$ be
the time spent by a particle performing a random walk generated by
matrix $A$ until the first hitting a point $z\in\mathbb{Z}^{d}$. In
a similar way, $\tau^{-}_{z}$ is the time spent by the particle {\it
after leaving the starting point} of the random walk until the first
hitting the point~$z$. If the starting point of the random walk is
$z$ then the first hitting $z$ means the first return to $z$. Denote
by $H_{x,z}(t)$ and $H^{-}_{x,z}(t)$, $t\geq0$, the (improper)
c.d.f. of $\tau_{z}$ and $\tau^{-}_{z}$, respectively, given that
the starting point of the random walk is $x\in\mathbb{Z}^{d}$.
Obviously, ${H_{x,z}(t)=G_{2}\ast H^{-}_{x,z}(t)}$ for $t\geq0$ and
$x,z\in\mathbb{Z}^{d}$. Set also ${K(t):=\alpha
f'(1)G_{1}(t)+(1-\alpha)G_{1}\ast H^{-}_{{\bf 0},{\bf 0}}(t)}$ and
$h(s):=\alpha(f(1-s)-1+f'(1)s)$, $s\in[0,1]$. Note that the function
$K(t)$ and the function $K_{d}(t)$, $d\in\mathbb{N}$, arising in
\cite{VT_4} and \cite{VT_Siberia}, coincide for each $t\geq0$. Thus,
Lemma 2.3 in \cite{VT_4} and Lemma 11 in \cite{VT_Siberia} in which
the asymptotic properties of c.d.f. $K_{d}(t)$ and its density
$k_{d}(t)$ are established, as $t\to\infty$, may be applied to our
function $K$.

\begin{Lm}\label{L:d=4}
For $y\in\mathbb{Z}^{d}$, $y\neq{\bf 0}$, one has
\begin{equation}\label{equation_d=4}
q(t;{\bf 0},y)=(1-\alpha)G_{1}\ast\left(H^{-}_{{\bf 0},y}(t)-H_{{\bf
0},y}(t)\right)+q(\cdot;{\bf 0},y)\ast K(t)-h(q(\cdot;{\bf
0},y))\ast G_{1}(t).
\end{equation}
\end{Lm}
{\sc Proof. }Recall integral equations (see, e.g., \cite{Sev}, Ch.8,
Sec.1) for probability generating functions
${\vec{F}(t;\vec{s}\,):=(F_{1}(t;\vec{s}\,),\ldots,F_{6}(t;\vec{s}\,))}$
of a six-dimensional Bellman-Harris process
$$F_{i}(t;\vec{s}\,)=s_{i}(1-G_{i}(t))+\int\nolimits_{0}^{t}{f_{i}\left(\vec{F}(t-u;\vec{s}\,)\right)
\,d G_{i}(u)},\quad t\geq0,\quad s_{i}\in[0,1],\quad i=1,\ldots,6.$$
By setting here $\vec{s}=(1,0,1,1,1,1)$ and substituting the
explicit formulae for $G_{j},$ $j=3,4,5,6$, and $f_{i}$,
$i=1,\ldots,6$, we get six integral equations in functions
$F_{i}(t):=F_{i}(t;(1,0,1,1,1,1))$, $t\geq0$, $i=1,\ldots,6$.
Substituting the fourth and the sixth ones into the second equation
and solving the obtained renewal equation in $F_{2}(\cdot)$ we find
\begin{eqnarray*}
F_{2}(t)&=&G_{2}\ast\left(1-H^{-}_{{\bf 0},y,{\bf 0}}(t)-H^{-}_{{\bf
0},{\bf 0},y}\ast\sum\nolimits_{k=0}^{\infty} {H^{\ast k}_{{\bf
0},{\bf 0},y}(t)}\right)\\
&+&F_{1}\ast G_{2}\ast\left(H^{-}_{{\bf 0},y,{\bf
0}}(0)+\left(H^{-}_{{\bf 0},y,{\bf 0}}(\cdot)-H^{-}_{{\bf 0},y,{\bf
0}}(0)\right)\ast\sum\nolimits_{k=0}^{\infty}{H^{\ast k}_{{\bf
0},{\bf 0},y}(t)}\right)
\end{eqnarray*}
where $H_{x,z,r}(t):=G_{2}\ast H^{-}_{x,z,r}(t)$, $t\geq0$,
$x,z,r\in\mathbb{Z}^{d}$, $z\neq r$. Now we substitute the last
equation as well as the third and the fifth equations in functions
$F_{i}$ into the first one. After some algebraic transformations we
obtain the following non-linear integral equation in function
$F_{1}$
\begin{eqnarray}
F_{1}(t)&=&1-\alpha
G_{1}\ast(1-f(F_{1}(t)))-(1-\alpha)G_{1}\ast(1-F_{1}(\cdot))\ast
H^{-}_{{\bf 0},{\bf 0}}(t)\nonumber\\
&-&(1-\alpha)G_{1}\ast(H^{-}_{{\bf 0},y}(t)-H_{{\bf
0},y}(t))\label{equation_d=4_F1(t)}
\end{eqnarray}
provided that the following two equalities are valid
$$H^{-}_{{\bf 0},{\bf 0}}(t)=H^{-}_{{\bf 0},{\bf 0},y}(t)+\sum\nolimits_{k=0}^{\infty}{H^{-}_{{\bf
0},y,{\bf 0}}\ast H^{\ast k}_{y,y,{\bf 0}}\ast H_{y,{\bf
0},y}(t)},\quad H^{-}_{{\bf 0},y}(t)=H^{-}_{{\bf 0},y,{\bf
0}}\ast\sum\nolimits_{k=0}^{\infty}{H^{\ast k}_{{\bf 0},{\bf
0},y}(t)},$$ for each $t\geq0$. The first of them is true since any
trajectory from ${\bf 0}$ to ${\bf 0}$ of a particle performing a
random walk on $\mathbb{Z}^{d}$ either passes $y$ exactly $k$ times,
$k=1,2,\ldots,$ or does not hit $y$ until the first returning to
${\bf 0}$. Similar argument justifies the second equality as well.
Recall that due to the connection between the CBRW on
$\mathbb{Z}^{d}$ and the constructed Bellman-Harris process one has
$q(t;{\bf 0},y)={\sf P}_{1}(Z_{2}(t)>0)=1-F_{1}(t)$. Hence,
rewriting \eqref{equation_d=4_F1(t)} as an equation in $q(t;{\bf
0},y)$ we come to \eqref{equation_d=4}. $\square$

\section{Proofs of Theorems \ref{T:q(t;x,y)sim} and \ref{T:limits}}

First of all, we derive some integral equations to be treated in
this section. Consider equation~\eqref{bKe_F(s,t;x,y)} as
inhomogeneous one for differential equation \eqref{bKe_m(t;x,y)} in
Banach space $l_{\infty}(\mathbb{Z}^{d})$. By the variation of
constant formula we infer (for a similar deduction see \cite{B_TSP})
that
\begin{equation}\label{Q(s,t;x,y)_integral equation}
q(s,t;x,y)=(1-s)m(t;x,y)-\int\nolimits_{0}^{t}{m(t-u;x,{\bf
0})h(q(s,u;{\bf 0},y))\,du}
\end{equation}
where $q(s,t;x,y)=1-F(s,t;x,y)$, $s\in[0,1]$, $t\geq0$,
$x,y\in\mathbb{Z}^{d}$. Substituting $x={\bf 0}$ in the last
equation we come to an integral equation in function $q(s,t;{\bf
0},y)$
\begin{equation}\label{Q(s,t;0,y)_integral equation}
q(s,t;{\bf 0},y)=(1-s)m(t;{\bf
0},y)-\int\nolimits_{0}^{t}{m(t-u;{\bf 0},{\bf 0})h(q(s,u;{\bf
0},y))\,du}.
\end{equation}
Note that $q(0,t;x,y)$ is equal to $q(t;x,y)$. Thus, on account of
\eqref{Q(s,t;x,y)_integral equation} one has
\begin{equation}\label{Q(t;x,y)_integral equation}
q(t;x,y)=m(t;x,y)-\int\nolimits_{0}^{t}{m(t-u;x,{\bf 0})h(q(u;{\bf
0},y))\,du}.
\end{equation}
Substituting $x={\bf 0}$ in \eqref{Q(t;x,y)_integral equation} we
derive an integral equation in function $q(t;{\bf 0},y)$
\begin{equation}\label{Q(t;0,y)_integral equation}
q(t;{\bf 0},y)=m(t;{\bf 0},y)-\int\nolimits_{0}^{t}{m(t-u;{\bf
0},{\bf 0})h(q(u;{\bf 0},y))\,du}.
\end{equation}

Now let us prove Theorems \ref{T:q(t;x,y)sim} and \ref{T:limits} for
$x={\bf 0}$. Since their proofs depend on $d\in\mathbb{N}$
essentially, we have to consider the cases $d=1$, $d=2$, $d=3$,
$d=4$ and $d\geq5$ separately. Evidently, Theorem
\ref{T:q(t;x,y)sim} for $x={\bf 0}$ and $d\geq6$ is implied by Lemma
\ref{L:d>=6}. Due to Lemmas
\ref{L:monotonicity_m(t;y,y)}--\ref{L:d=5} and
equation~\eqref{Q(t;0,y)_integral equation} the proof of Theorem
\ref{T:q(t;x,y)sim} for $x={\bf 0}$ in the respective cases $d=1$,
$d=2$, $d=3$ and $d=5$ mainly follows the scheme proving,
respectively, Theorem 2 in \cite{VTY}, Theorem 2 in \cite{B_TVP},
Theorem 4 in \cite{VT_Siberia} (item 3) and Theorem 4 in
\cite{VT_Siberia} (item 4). Moreover, by virtue of Lemma \ref{L:d=4}
the proof of Theorem \ref{T:q(t;x,y)sim} for $x={\bf 0}$ and $d=4$
is similar to that of Theorem 1.1 in \cite{VT_4}. So, we give only a
few comments on the proof of Theorem \ref{T:q(t;x,y)sim} for $x={\bf
0}$ and $d\leq5$.

If $d=1$ then the equality $\int\nolimits_{0}^{\infty}{h(q(u;{\bf
0},y))\,du}=(1-\alpha)a^{-1}$ is valid. Furthermore, in view of
\eqref{p(t;x,y)sim,p(t;0,0)-p(t;x,y)}, \eqref{int_bKe_m(t;x,y)},
\eqref{int_fKe_m(t;x,y)} and Theorem 5 in \cite{VT_Siberia} one gets
the useful estimate
$$m(t;{\bf 0},y)-(1-\alpha)a^{-1} m(t;{\bf 0},{\bf
0})=O\left(t^{-3/2}\right),\quad t\to\infty.$$ When $d=2$ one can
check the strict inequality
$J(0;y)=\int\nolimits_{0}^{\infty}{h(q(u;{\bf
0},y))\,du}<(1-\alpha)a^{-1}$. However, if $d=3$ then
$\int\nolimits_{0}^{\infty}{h(q(u;{\bf
0},y))\,du}=(1-\alpha)a^{-1}G_{0}({\bf 0},y)G^{-1}_{0}({\bf 0},{\bf
0})$ and
$$m(t;{\bf 0},y)-(1-\alpha)a^{-1}G_{0}({\bf 0},y)G^{-1}_{0}({\bf 0},{\bf
0})\,m(t;{\bf 0},{\bf 0})=O(t^{-1}),\quad t\to\infty,$$ by virtue of
\eqref{p(t;x,y)sim,p(t;0,0)-p(t;x,y)}, \eqref{int_fKe_m(t;x,y)} and
Theorem 5 in \cite{VT_Siberia}. For $d=4$ the first summand in
\eqref{equation_d=4} is $o(t^{-1})$, $t\to\infty$, by Lemma 3 in
\cite{LMJ} and it does not contribute to the (main term of)
asymptotic behavior of $q(t;{\bf 0},y)$. As for $d=5$, one has
$\int\nolimits_{0}^{\infty}{h(q(u;{\bf
0},y))\,du}=(1-\alpha)a^{-1}G_{0}({\bf 0},y)G^{-1}_{0}({\bf 0},{\bf
0})$ and
$$m(t;{\bf
0},y)-(1-\alpha)a^{-1}G_{0}({\bf 0},y)G^{-1}_{0}({\bf 0},{\bf
0})\,m(t;{\bf 0},{\bf 0})=O\left(t^{-3/2}\right),\quad t\to\infty,$$
in view of \eqref{p(t;x,y)sim,p(t;0,0)-p(t;x,y)},
\eqref{int_fKe_m(t;x,y)}, Theorem 5 and Corollary 1 in
\cite{VT_Siberia}. Thus, Theorem \ref{T:q(t;x,y)sim} is proved for
$x={\bf 0}$.

Turn to Theorem \ref{T:limits} when $x={\bf 0}$. The proof of
Theorem \ref{T:limits} for $x={\bf 0}$ is similar to those of
Theorem 4 in \cite{VT_DM}, Theorem 2 in \cite{B_TVP} and Theorem 4
in \cite{LMJ} for $d=1,3$, $d=2$ and $d\geq5$, respectively. Note
only that the constant $c^{\ast}$ arising in the proof of Theorem
\ref{T:limits} for $x={\bf 0}$ in contrast to its counterpart in
Theorem 4 in \cite{VT_DM} is equal to
$\sigma^{2}\gamma^{2}_{1}a/(2(1-\alpha))$ and
$\sigma^{2}aG^{3}_{0}({\bf 0},{\bf 0})G_{0}({\bf
0},y)/(8\pi^{2}\gamma^{2}_{3}(1-\alpha))$ when $d=1$ and $d=3$,
respectively. At last, the constant $c^{\ast}_{d}$ appearing in
Theorem 4 in \cite{LMJ} equals $(1-\alpha)G_{0}({\bf
0},y)\sigma^{2}/\left(2a\,G_{0}({\bf 0},{\bf 0})m^{2}_{d}\right)$ in
the case of Theorem \ref{T:limits} for $x={\bf 0}$ and $d\geq5$.
Since the limit theorem for $\mu(t;{\bf 0})$ when $d=4$ was
established by another approach, namely the moment method, we give
the detailed proof of the limit theorem for $\mu(t;y)$ when $d=4$.
So, to complete the proof of Theorem \ref{T:limits} for $x={\bf 0}$
we dwell on the case $d=4$ in detail.

Set $s(t):=s(t;\lambda)=\exp\{-\lambda\ln^{2}{t}/(c^{\ast} t)\}$
where $c^{\ast}:=\sigma^{2} a G^{3}_{0}({\bf 0},{\bf 0}) G_{0}({\bf
0},y)/(3\gamma^{2}_{4}(1-\alpha))$, $t>0$ and $\lambda\geq0$. By
Theorems \ref{T:m(t;x,y)} and \ref{T:q(t;x,y)sim} for $x={\bf 0}$
and $d=4$ we see that
\begin{equation}\label{E(mu(y,y)|mu(t,y)>0)sim}
{\sf E}_{\bf 0}\left(\left.\mu(t;y)\right|\mu(t;y)>0\right)=
\frac{m(t;{\bf 0},y)}{q(t;{\bf 0},y)}\sim\frac{c^{\ast}
t}{\ln^{2}{t}},\quad t\to\infty.
\end{equation}
The inequality $1-e^{-z}\leq z$ for $z\geq0$ yields
$$q(s(t),u;{\bf 0},y)={\sf E}_{\bf
0}\left(1-\exp\left\{-\frac{\lambda\ln^{2}{t}\,\mu(u;y)}{c^{\ast}t}\right\}\right)
\leq\frac{\lambda\ln^{2}{t}}{c^{\ast}t}{\sf E}_{\bf
0}\mu(u;y)=\frac{\lambda\ln^{2}{t}}{c^{\ast}t}m(u;{\bf 0},y)$$ where
$u\geq0$ and $t>0$. By virtue of this estimate combined with
Theorem~1 and the inequality $h(z)\leq\sigma^{2}z^{2}$ (being true
for $z\geq0$ small enough) one has for $t$ large enough
\begin{eqnarray}
& &\int\nolimits_{0}^{t/\ln^{3}{t}}{m(t-u;{\bf 0},{\bf
0})h(q(s(t),u;{\bf 0},y))du}\nonumber\\
&\leq&\frac{\sigma^{2}\lambda^{2}\ln^{4}{t}}{c^{\ast
2}t^{2}}\int\nolimits_{0}^{t/\ln^{3}{t}}{m^{2}(u;{\bf
0},y)m(t-u;{\bf 0},{\bf
0})du}=\frac{\rho_{1}(t;\lambda)\ln{t}}{t}.\label{estimates_of_tails_1}
\end{eqnarray}
Here $\rho_{1}\in\mathcal{U}$ and $\mathcal{U}$ is the class of all
bounded functions $\rho(t;\lambda)$ vanishing as $t\to\infty$
uniformly in $\lambda\in[0,b]$, whatever positive $b$ is taken. In a
similar way, we obtain
\begin{eqnarray}
& &\int\nolimits_{t-t/\ln^{2}{t}}^{t}{m(t-u;{\bf 0},{\bf
0})h(q(s(t),u;{\bf
0},y))du}\nonumber\\
&\leq&\frac{\sigma^{2}\lambda^{2}\ln^{4}{t}}{c^{\ast
2}t^{2}}\int\nolimits_{t-t/\ln^{2}{t}}^{t}{m^{2}(u;{\bf
0},y)m(t-u;{\bf 0},{\bf
0})du}=\frac{\rho_{2}(t;\lambda)\ln{t}}{t}\label{estimates_of_tails_2}
\end{eqnarray}
for $\rho_{2}\in\mathcal{U}$. It is not difficult to show that
uniformly in $u\in[t/\ln^{3}{t},t-t/\ln^{2}{t}]$
\begin{equation}\label{useful_relations}
\ln{u}\sim\ln{t},\quad\ln(t-u)\sim\ln{t},\quad t\to\infty.
\end{equation}
These facts, Theorem \ref{T:m(t;x,y)} and the relation
$h(z)\sim\sigma^{2}z^{2}/2,\;z\to0,$ allow us to claim that
\begin{eqnarray*}
I(t;\lambda):&=&\int\nolimits_{t/\ln^{3}{t}}^{t-t/\ln^{2}{t}}{h(q(s(t),u;{\bf
0},y))m(t-u;{\bf 0},{\bf 0})du}\\
&=&\frac{\sigma^{2}m(t;{\bf 0},{\bf
0})}{2}\int\nolimits_{t/\ln^{3}{t}}^{t-t/\ln^{2}{t}}{q^{2}(s(t),u;{\bf
0},y)du}\;(1+\rho_{3}(t;\lambda))
\end{eqnarray*}
where $\rho_{3}\in\mathcal{U}$. After changing the variable
$u=t\upsilon$ and using Theorems 1 and 2 for $x={\bf 0}$ and $d=4$
we get
\begin{equation}\label{I(t)2_d=4}
I(t;\lambda)=\frac{3}{2q(t;{\bf\,0},y)}\int\nolimits_{1/\ln^{3}{t}}^{1-1/\ln^{2}{t}}{q^{2}(s(t),t\upsilon;{\bf
0},y)\,d\upsilon}\;(1+\rho_{4}(t;\lambda)),\quad\rho_{4}\in\mathcal{U}.
\end{equation}
In the last integral the function $q(s(t;\lambda),t\upsilon;{\bf
0},y)$ can be replaced by
$q(s(t\upsilon;\lambda\upsilon),t\upsilon;{\bf 0},y)$. Indeed, as
$1-e^{-z}\leq z$ for $z\geq0$, we have
\begin{eqnarray*}
& &|q(s(t;\lambda),t\upsilon;{\bf
0},y)-q(s(t\upsilon;\lambda\upsilon),t\upsilon;{\bf 0},y)|\\
&=&{\sf E}_{\bf
0}\left(\exp\left\{-\frac{\lambda\upsilon\ln^{2}(t\upsilon)}{c^{\ast}t\upsilon}\mu(t\upsilon;y)\right\}
-\exp\left\{-\frac{\lambda\ln^{2}{t}}{c^{\ast}t}\mu(t\upsilon;y)\right\}\right)\\
&\leq&{\sf E}_{\bf
0}{\left(1-\exp\left\{-\frac{\lambda(-2\ln{t}\ln{\upsilon}-\ln^{2}{\upsilon})}
{c^{\ast}t}\mu(t\upsilon;y)\right\}\right)}\leq\frac{\lambda(-2\ln{t}\ln{\upsilon}-\ln^{2}{\upsilon})}
{c^{\ast}t}m(t\upsilon;{\bf 0},y).
\end{eqnarray*}
Since functions $z\ln{z}$ and $z\ln^{2}{z}$ are bounded for
$z\in(0,1)$, by virtue of Theorems \ref{T:m(t;x,y)} and
\ref{T:q(t;x,y)sim} for $x={\bf 0}$ along with relation
\eqref{useful_relations} we see that uniformly in
$\upsilon\in[1/\ln^{3}{t},1-1/\ln^{2}{t}]$ and
$0\leq\lambda\leq\Lambda$ with an arbitrary positive $\Lambda$
\begin{equation}\label{frac-frac=o(1)}
\frac{q(s(t;\lambda),t\upsilon;{\bf 0},y)}{q(t\upsilon;{\bf
0},y)}-\frac{q(s(t\upsilon;\lambda\upsilon),t\upsilon;{\bf
0},y)}{q(t\upsilon;{\bf 0},y)}\to0,\quad t\to\infty.
\end{equation}
Set $\varphi(t;\lambda):=q(s(t;\lambda),t;{\bf 0},y)/(\lambda
q(t;{\bf 0},y))$, $t>0$, $\lambda\geq0$. Then dividing both sides of
\eqref{Q(s,t;0,y)_integral equation} by $\lambda q(t;{\bf 0},y)$ and
using \eqref{estimates_of_tails_1}--\eqref{frac-frac=o(1)} along
with Theorem \ref{T:q(t;x,y)sim} for $x=0$ and relation
$1-e^{-z}\sim z$, $z\to0$, we obtain
$$\varphi(t;\lambda)=1+\rho_{5}(t;\lambda)-\frac{3\lambda}{2}
\int\nolimits_{1/\ln^{3}{t}}^{1-1/\ln^{2}{t}}{\varphi^{2}(t\upsilon;\lambda\upsilon)
\,d\upsilon},\quad
\rho_{5}\in\mathcal{U}.$$ Changing the variable $w=\lambda\upsilon$
leads to the following relation
$$\varphi(t;\lambda)=1+\rho_{5}(t;\lambda)-\frac{3}{2}
\int\nolimits_{\lambda/\ln^{3}{t}}^{\lambda(1-1/\ln^{2}{t})}{\varphi^{2}
\left(\frac{t\,w}{\lambda};w\right)\,d w}.$$ The argument similar to
the proof of Theorem 4 in \cite{VT_DM} establishes that
\begin{equation}\label{varphi_to_varphi}
\lim\limits_{t\to\infty}{\varphi(t;\lambda)}=
\varphi(\lambda)=\frac{2}{3\lambda+2},\quad0<\lambda\leq\Lambda_{0},
\end{equation}
where $\Lambda_{0}$ is some positive number and $\varphi(\lambda)$
is the unique solution of the equation
$$\varphi(\lambda)
=1-\frac{3}{2}\int\nolimits_{0}^{\lambda}{\varphi^{2}(w)\,d
w},\quad\lambda\geq0.$$ Invoking the definition of
$\varphi(t;\lambda)$ we rewrite relation \eqref{varphi_to_varphi} by
way of
\begin{equation}\label{LT}
\lim\limits_{t\to\infty}{{\sf E}_{\bf
0}\left\{\left.\exp\left\{-\frac{\lambda\ln^{2}{t}\,\mu(t;y)}{c^{\ast}t}\right\}
\right|\mu(t;y)>0\right\}}=
1-\lambda\lim\limits_{t\to\infty}{\varphi(t;\lambda)}=\frac{1}{3}+\frac{2}{3}
\cdot\frac{2}{3\lambda+2}
\end{equation}
for $0<\lambda\leq\Lambda_{0}$. Since both the Laplace transform of
a non-negative random variable and the function
$1/3+2/3\cdot2/(3\lambda+2)$ are analytic and bounded in the domain
$\{\lambda: Re\,\lambda>0\}\subset\mathbb{C}$, by the uniqueness
theorem for analytic functions relation \eqref{LT} is valid for each
$\lambda$ with $Re\,\lambda>0$ (for an analogous deduction see,
e.g., \cite{V_1986}). Combining \eqref{E(mu(y,y)|mu(t,y)>0)sim} and
\eqref{LT} we complete the proof of Theorem \ref{T:limits} for $x=0$
and $d=4$. Thus, Theorem \ref{T:limits} is proved for $x={\bf 0}$.

Next we prove Theorems \ref{T:q(t;x,y)sim} and \ref{T:limits} when
$x\neq{\bf 0}$. As a preliminary we derive some more integral
equations. In the framework of CBRW on $\mathbb{Z}^{d}$, the parent
particle can either hit the point ${\bf 0}$ or not within time
interval $[0,t]$. In the latter case at time $t$ there is a single
particle on $\mathbb{Z}^{d}$ located at the point $y$ or outside it.
Consequently,
\begin{eqnarray}
{\sf E}_{x}{s^{\mu(t;y)}}&=&{\sf
E}_{x}{s^{\mu(t;y)}\mathbb{I}(\tau_{{\bf 0}}\leq t)}+{\sf
E}_{x}{s^{\mu(t;y)}\mathbb{I}(\tau_{{\bf 0}}>t,\mu(t;y)=1)}+{\sf
E}_{x}{s^{\mu(t;y)}\mathbb{I}(\tau_{{\bf 0}}>t,\mu(t;y)=0)}\nonumber\\
&=&{\sf E}_{x}{s^{\mu(t;y)}\mathbb{I}(\tau_{{\bf 0}}\leq t)}+s{\sf
P}_{x}(\tau_{{\bf 0}}>t,\mu(t;y)=1)+{\sf P}_{x}(\tau_{{\bf
0}}>t,\mu(t;y)=0)\label{new_equation_x}
\end{eqnarray}
where $\mathbb{I}(\cdot)$ stands for the indicator of a set.
Evidently, the first summand in \eqref{new_equation_x} can be
rewritten in the form
\begin{eqnarray}\label{new_equation_x_1_summand}
{\sf E}_{x}{s^{\mu(t;y)}\mathbb{I}(\tau_{{\bf 0}}\leq
t)}&=&\int\limits_{\{\tau_{{\bf 0}}\leq t\}}{s^{\mu(t;y)}\,d\,{\sf
P}_{x}}=\int\limits_{\{\tau_{{\bf 0}}\leq t\}}{{\sf
E}_{x}\left(\left.s^{\mu(t;y)}\right|\tau_{{\bf 0}}\right)\,d\,{\sf
P}_{x}}\nonumber\\
&=&\int\nolimits_{0}^{t}{{\sf
E}_{x}\left(\left.s^{\mu(t;y)}\right|\tau_{{\bf 0}}=u\right)\,d
H_{x,{\bf 0}}(u)}=\int\nolimits_{0}^{t}{{\sf E}_{{\bf
0}}{s^{\mu(t-u;y)}}\,d H_{x,{\bf 0}}(u)}.
\end{eqnarray}
It is easily seen that the probability at the second summand in
\eqref{new_equation_x} can be represented as follows
\begin{eqnarray}
{\sf P}_{x}\left(\tau_{{\bf 0}}>t,\mu(t;y)=1\right)&=&H_{x,y,{\bf
0}}\ast\sum\nolimits_{k=0}^{\infty}{H^{\ast k}_{y,y,{\bf
0}}}\ast(1-G_{2}(t))\quad\mbox{when}\quad x\neq y,\label{new_equation_x_2_summand}\\
{\sf P}_{y}\left(\tau_{{\bf
0}}>t,\mu(t;y)=1\right)&=&\sum\nolimits_{k=0}^{\infty}{H^{\ast
k}_{y,y,{\bf
0}}}\ast(1-G_{2}(t)).\label{new_equation_x_2_summand_x=y}
\end{eqnarray}
It also turns convenient to write the third summand in
\eqref{new_equation_x} in the form
\begin{eqnarray}\label{new_equation_x_3_summand}
{\sf P}_{x}(\tau_{{\bf 0}}>t,\mu(t;y)=0)&=&1-H_{x,{\bf
0}}(t)-H_{x,y,{\bf 0}}\ast\sum\nolimits_{k=0}^{\infty}{H^{\ast
k}_{y,y,{\bf 0}}}\ast(1-G_{2}(t))\quad\mbox{if}\quad x\neq y,\quad\label{new_equation_x_3_summand}\\
{\sf P}_{y}(\tau_{{\bf 0}}>t,\mu(t;y)=0)&=&1-H_{y,{\bf
0}}(t)-\sum\nolimits_{k=0}^{\infty}{H^{\ast k}_{y,y,{\bf
0}}}\ast(1-G_{2}(t)).\label{new_equation_x_3_summand_x=y}
\end{eqnarray}
Combining relations
\eqref{new_equation_x}--\eqref{new_equation_x_3_summand_x=y} we come
to the desired integral equations
\begin{eqnarray*}
q(s,t;x,y)\!\!\!&=&\!\!\!(1-s)H_{x,y,{\bf
0}}\ast\sum\nolimits_{k=0}^{\infty}{H^{\ast k}_{y,y,{\bf
0}}}\ast(1-G_{2}(t))+\int\nolimits_{0}^{t}{q(s,t-u;{\bf 0},y)\,d
H_{x,{\bf
0}}(u)}\;\;\mbox{if}\; x\neq y,\\
q(s,t;y,y)\!\!\!&=&\!\!\!(1-s)\sum\nolimits_{k=0}^{\infty}{H^{\ast
k}_{y,y,{\bf 0}}}\ast(1-G_{2}(t))+\int\nolimits_{0}^{t}{q(s,t-u;{\bf
0},y)\,d H_{y,{\bf 0}}(u)}.
\end{eqnarray*}
In particular, for $s=0$ one has
\begin{eqnarray}
q(t;x,y)\!\!\!&=&\!\!\!H_{x,y,{\bf
0}}\ast\sum\nolimits_{k=0}^{\infty}{H^{\ast k}_{y,y,{\bf
0}}}\ast(1-G_{2}(t))+\int\nolimits_{0}^{t}{q(t-u;{\bf 0},y)\,d
H_{x,{\bf
0}}(u)}\;\;\mbox{when}\; x\neq y,\quad\label{start_any_x_s=0}\\
q(t;y,y)\!\!\!&=&\!\!\!\sum\nolimits_{k=0}^{\infty}{H^{\ast
k}_{y,y,{\bf 0}}}\ast(1-G_{2}(t))+\int\nolimits_{0}^{t}{q(t-u;{\bf
0},y)\,d H_{y,{\bf 0}}(u)}.\quad\label{start_any_x_x=y_s=0}
\end{eqnarray}

Now we have the tools for proving Theorems \ref{T:q(t;x,y)sim} and
\ref{T:limits} for $x\neq{\bf 0}$. To establish Theorem
\ref{T:q(t;x,y)sim} for $x\neq{\bf 0}$ and $d\neq2$ we employ
equations \eqref{start_any_x_s=0} and \eqref{start_any_x_x=y_s=0}.
It is not difficult to see that the first summands in the right side
of \eqref{start_any_x_s=0} and \eqref{start_any_x_x=y_s=0} are equal
to $p(t;x,y)-\int\nolimits_{0}^{t}{p(t-u;{\bf 0},y)\,d H_{x,{\bf
0}}(u)}$ for $x\neq y$ and $x=y$, respectively. The latter
expression can be rewritten as follows
\begin{eqnarray}
& &p(t;x,y)-\int\nolimits_{0}^{t}{p(t-u;{\bf 0},y)d H_{x,{\bf
0}}(u)}\nonumber\\
&=&p(t;x,y)-p(t;x,{\bf 0})+\int\nolimits_{0}^{t}{(p(t-u;{\bf 0},{\bf
0})-p(t-u;{\bf 0},y))d H_{x,{\bf 0}}(u)}\label{x_neq_0}
\end{eqnarray}
due to the obvious relation $p(t;x,{\bf
0})=\int\nolimits_{0}^{t}{p(t-u;{\bf 0},{\bf 0})\,d H_{x,{\bf
0}}(u)}$. The asymptotic behavior of the first summand at the
right-hand side of \eqref{x_neq_0} is given by formula
\eqref{p(t;x,y)sim,p(t;0,0)-p(t;x,y)} whereas the asymptotic
behavior of the second summand in \eqref{x_neq_0} can be found with
the help of relation \eqref{p(t;x,y)sim,p(t;0,0)-p(t;x,y)}, Lemma 3
in \cite{LMJ} and Lemma 5.1.2 in \cite{Y_Book}. Finally, the first
summands in \eqref{start_any_x_s=0} and \eqref{start_any_x_x=y_s=0}
are $O(t^{-3/2})$ when $d=1$ and $O(t^{-d/2})$ when $d\geq3$. Hence,
the first summands in \eqref{start_any_x_s=0} and
\eqref{start_any_x_x=y_s=0} are $o(q(t;{\bf 0},y))$, as
$t\to\infty$, by Theorem \ref{T:q(t;x,y)sim} for $x={\bf 0}$.
Moreover, on account of Lemma 3 in \cite{LMJ} and Lemma~5.1.2 in
\cite{Y_Book} we reveal that the last summands in
\eqref{start_any_x_s=0} and \eqref{start_any_x_x=y_s=0} are
equivalent to $q(t;{\bf 0},y)$ and $q(t;{\bf 0},y)\,G_{0}(x,{\bf
0})G^{-1}_{0}({\bf 0},{\bf 0})$ for $d=1$ and $d\geq3$,
respectively, as $t\to\infty$. Hence Theorem \ref{T:q(t;x,y)sim} is
proved for $x\neq{\bf 0}$ and $d\neq2$. As for Theorem
\ref{T:q(t;x,y)sim} when $x\neq{\bf 0}$ and $d=2$ as well as Theorem
\ref{T:limits} for $x\neq{\bf 0}$, we only note that their proofs
bear on analysis of equations \eqref{Q(s,t;x,y)_integral equation}
and \eqref{Q(t;x,y)_integral equation}. Since the proofs are similar
to that of Theorem 5 in \cite{LMJ}, they are omitted. So, Theorems
\ref{T:q(t;x,y)sim} and \ref{T:limits} are proved completely.
$\square$

\vskip0.5cm The author is grateful to Associate Professor
E.B.Yarovaya for permanent attention and to Professor V.A.Vatutin
for valuable remarks. Special thanks are to Professors I.Kourkova
and G.Pag\`es for invitation to LPMA UPMC (Paris-6) where this work
was started.

\newpage
Ekaterina Vl. BULINSKAYA, \vskip0,2cm Faculty of Mathematics and
Mechanics,

Lomonosov Moscow State University,

Moscow 119991, Russia

\vskip0,5cm

{\it E-mail address}: bulinskaya@yandex.ru

\end{document}